# A New Way to Proof 3x+1 Problem

Yanlong Zhou[1]

1. Beijing Inst Technol, Dept Transportat Engn, Beijing 100081, Peoples R China

ABSTRACT：Under the 3x+1 problem, classified the number into four kind by mod 4.

$$z=\begin{cases}\alpha = 4*k+1\\ \beta = 4*k+2\\ \eta = 4*k+3\\ \gamma = 4*k+4\end{cases} \quad k=0,1,2,\cdots$$

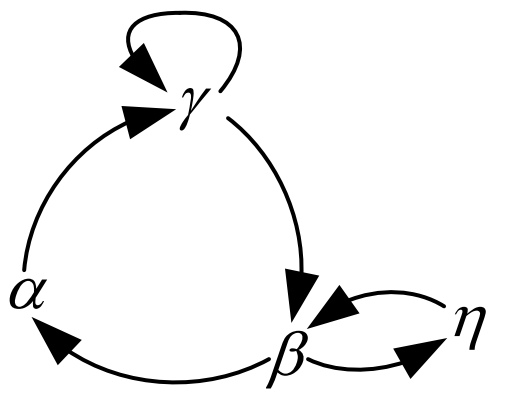

The four kind number can form a cycle base on 3x+b1 problem. Base on this cycle, if the number of $\eta$ kind number is zero the 3x+1 will be proofed.

## 0 Introduction

Analysis for the numbers by the Collatz Problem

$$C(z)=\begin{cases}\dfrac{z}{2} & when \quad z=0(\mathrm{mod}\,2)\\ 3*z+1 & when \quad z=1(\mathrm{mod}\,2)\end{cases}$$

As we know the rule of Collatz Problem can express as $T(z)$

$$T(z)=\begin{cases}\dfrac{z}{2} & when \quad x=0(\mathrm{mod}\,2)\\ \dfrac{3*z+1}{2} & when \quad x=1(\mathrm{mod}\,2)\end{cases}$$

The 3x +1 problem is generally attributed to Lothar Collatz and has reportedly circulated since at least the early 1950's. Peter Eisele and Karl-Peter Hadeler (1990) studies iteration of the mappings f(x) = a + [x/b] on Z where a, b are positive integers[1]. Gaston Gonnet (1991) describes how to write computer code to efficiently compute 3x +1 function iterates for very large x using MAPLE[2].

Edward Belaga and Maurice Mignotte (2006) gives detailed tables of primitive cycles on the positive integers for 3x+d maps for 1 < d < 20000[3].Those great people do great amount work on the proof by the Glide Record, Delay Record, Residue, Completeness Record ,and Class Records.

But those just obey this rule, but if we obey this back rule, that mean we can get the all $N^+$ by '1'.

## 1 Analysis for the $N^+$ by back rule

If we summed that Collatz Problem is founded, just find the back way for the numbers, we can know that the '1' can establish the hold integers by the hold integers.

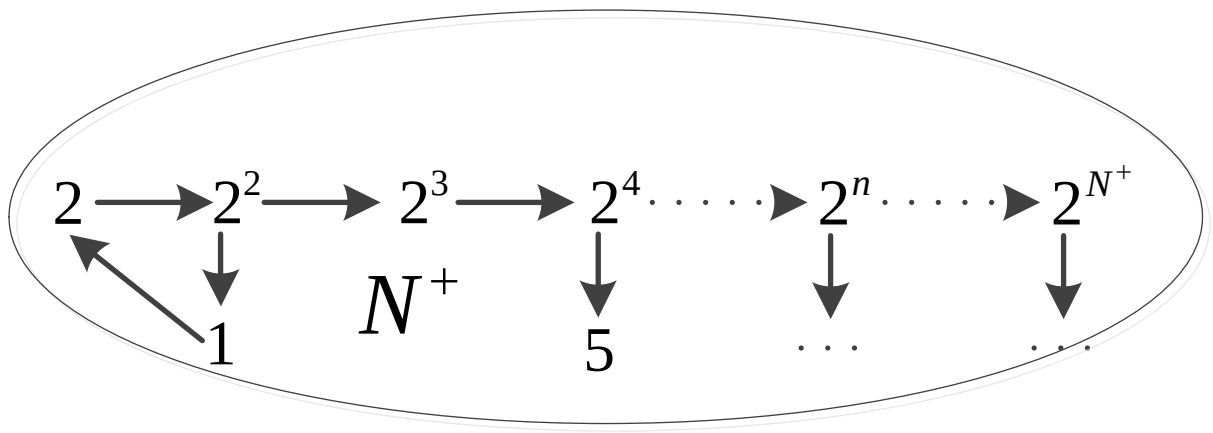

Fig1 Back Collatz Problem rule number's change path

The back for the rule by Collatz Problem is the ( $N^+$ , $N^+$ ,$\cdots$) map to $N^+$ , based on the fig 1. Because this number change rule by Collatz Problem based on the judgment number's parity, so rethinking the decimal numbers system generates is an important way for proofing Collatz Problem. When I learn the fractal. I find the fractal can give me a better view for the establishment of decimal numbers system.

## 2 Classify number into four kind

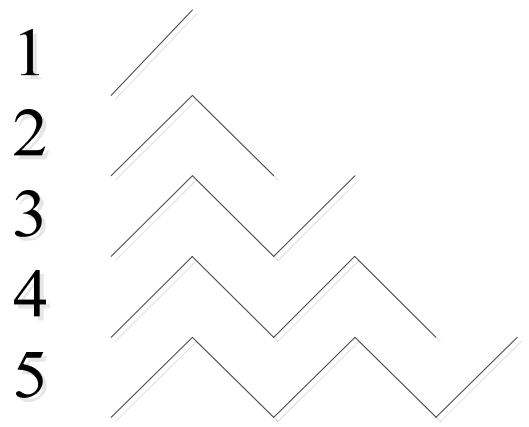

Fig.2 Polyline represents a number

If we classify the $N^+$ by the way for establishment of decimal numbers system, we can find four kinds of numbers by judging parity for the number of upper vertices ( $s$ ) and the number of under vertices ( $x$ ): the $\alpha$ kind number, $s_\alpha$ is an odd number, $x_\alpha$ is an odd number; the $\beta$ kind number, $s_\beta$ is an odd number, $x_\beta$ is an even number; the $\gamma$ kind number, $s_\gamma$ is an even number, $x_\gamma$ is an odd number; the $\eta$ kind number, $s_\eta$ is an even number, $x_\eta$ is an even number. So a decimal numbers system number z will be $z = x + s - 1$. An odd number $z_{odd}$ have $x_{odd} = s_{odd}$ , an even number $z_{even}$ have $x_{even} = s_{even} + 1$.

## 3 Build the cycle

And we can redefined the four kind number $\alpha, \beta, \gamma, \eta$ by another way:

$$z = \begin{cases} \alpha = 4*k+1 \\ \beta = 4*k+2 \\ \eta = 4*k+3 \\ \gamma = 4*k+4 \end{cases} \quad k = 0,1,2,\cdots \tag{1}$$

Use the $C(\mathrm{z})$, the will be

$$C(\mathrm{z})=\begin{cases} C(\alpha)=3*4*k+4 \\ C(\beta)=2*k+1 \\ C(\eta)=3*4*k+10 \\ C(\gamma)=2*k+2 \end{cases} \quad k=0,1,2,\cdots$$

And we set k be

$$k=\begin{cases} 2*l+1 \\ 2*l \end{cases} \quad l=0,1,2\cdots$$

Then we get

$$C(z)=\begin{cases} C(\alpha)=\begin{cases} 4*(3*2*l+3)+4=\gamma & k=2*l+1 \\ 3*4*(3*2*l)+4=\gamma & k=2*l \end{cases} \\ C(\beta)=\begin{cases} 4*l+3=\eta & k=2*l+1 \\ 4*l+1=\alpha & k=2*l \end{cases} \\ C(\eta)=\begin{cases} 4*(3*2*l+5)+2=\beta & k=2*l+1 \\ 4*(3*2*l+2)+2=\beta & k=2*l \end{cases} \\ C(\gamma)=\begin{cases} 4*l+4=\gamma & k=2*l+1 \\ 4*l+2=\beta & k=2*l \end{cases} \end{cases} \quad l=0,1,2,\cdots\cdots \tag{2}$$

As the Ep.6 show. We can get a graph for changing of the four kind number base on it.

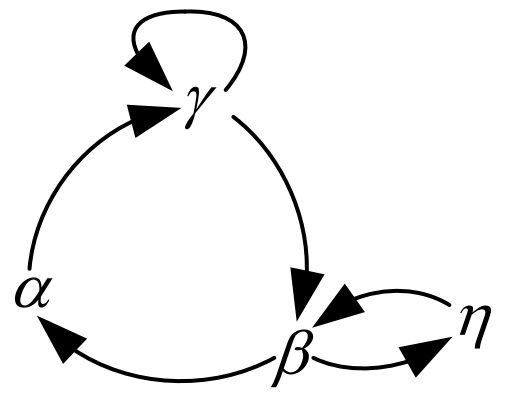

Fig.3 The four kind number cycle

So the 3z+1 problem will be the cycle system. The main cycle is $\alpha\to\gamma\to\beta\to\alpha$, the self-cycle $\gamma\to\gamma$ and two-cycle system $\beta\to\eta\to\beta$. if we want proof the 3z+1 problem.in an infinity time the main cycle is $1\to4\to2\to1$.so all two-cycle system is infinity , and self-cycle $\gamma\to\gamma$ is infinity too. It is obviously that the b self-cycle $\gamma\to\gamma$ is infinity. So we proof the two-cycle system $\beta\to\eta\to\beta$ is infinity.

## 4 Proof for two-cycle system $\beta\to\eta\to\beta$ is infinity

$$C(z)=\begin{cases} C(\beta)=\begin{cases}4*l+3=\eta & k=2*l+1\\ 4*l+1=\alpha & k=2*l\end{cases} \\ C(\eta)=\begin{cases}4*(3*2*l+5)+2=\beta & k=2*l+1\\ 4*(3*2*l+2)+2=\beta & k=2*l\end{cases}\end{cases} \quad l=0,1,2,\cdots\cdots$$

As we know the two-cycle system should be an infinity cycle, in another words, after infinity time we can get an end.

$$C^n(\beta)=4*h+1 \quad h=0,1,2,\cdots\cdots \quad n=1,2,3,\cdots\cdots \tag{3}$$

But if $C^n(\beta)=\alpha$, the trajectories n is an odd number.

$$\underbrace{\overbrace{\beta\to\eta\to}^{2}\underbrace{\beta\to\eta\to}_{2}\overbrace{\beta\to\alpha}^{1}}_{n}$$

Fig.4 The two-cycle system

So we let $n=2*m+1 \quad (m=0,1,2,\cdots\cdots)$, and the Ep.7 will be

$$4*h+1=C^{2*m+1}(\beta)=(\frac{\beta}{2})*(\frac{3}{2})^m+\sum_{i=1}^{m}\frac{3^{i-1}}{2^i} \quad m=0,1,2,\cdots\cdots$$

(4)

Simplify

$$\left(k+1\right)*3^m=(2*h+1)*2^m \quad m=0,1,2,\cdots\cdots \tag{5}$$

When $k$ is even, this function must be get a pier infinity number $(m,h)$ for $k$, because the left is odd, so that the right is odd too. Then we can get $m=0$, the right will be $2*h+1$, we can get a number $h$ for that. When $k$ is odd, we can let $k=2*t-1 \quad (t=1,2,3,\cdots\cdots)$

$$t*3^m=(2*h+1)*2^{m-1} \quad m=0,1,2,\cdots$$

So when $t$ is odd, the right must be odd too so $m=1$, then the right will be $2*h+1$ too, then we can find a number $h$ for that. When $t$ is even, we can let $t=2*w \quad (w=1,2,3,\cdots\cdots)$

$$w*3^m=(2*h+1)*2^{m-2} \quad m=0,1,2,\cdots$$

So when $w$ is odd, the right must be odd too so $m=2$, then the right will be $2*h+1$ too, then we can find a number $h$ for that. So we can know that every $k$ will have a pair infinity number $(m,h)$ for it. So we proof the two-cycle system is infinity.

## 5 Proof the main cycle is only by 1, 4, and 2.

We can know the finally of the main cycle $\alpha \to \gamma \to \beta \to \alpha$ based on 1, 4, and 2.

$$\begin{cases} \alpha = 4*h+1 \\ \beta = 4*k+2 \quad h,k,g = 0,1,2,\cdots\cdots \\ \gamma = 4*g+4 \end{cases}$$

$$\beta \to \alpha \quad (k+1)*3^m = (2*h+1)*2^m \quad m = 0,1,2,\cdots\cdots \tag{6}$$

$$\alpha \to \gamma \quad 4*g+4 = (4*h+1)*3+1 \tag{7}$$

$$\gamma \to \beta \quad 4*k+2 = \frac{4*g+4}{2^e} \quad e = 1,2,3,\cdots\cdots \tag{8}$$

As we know $4*h+1$ is $\alpha$, it will be $\beta$ .as the cycle $\alpha^0 \to \gamma^0 \to \beta^0 \to \alpha^1 \to \gamma^1 \to \beta^1$

$$4*h^0+1 = (\frac{4*k^0+2}{2})*(\frac{3}{2})^{m^0} + (\frac{3}{2})^{m^0} - 1 \quad m^0 = 0,1,2,\cdots\cdots$$

$$4*k^1+2 = \frac{(4*h^0+1)*3+1}{2^{e^0}} \quad e^0 = 1,2,3,\cdots\cdots$$

$$4*k^1+2 = \frac{((\frac{4*k^0+2}{2})*(\frac{3}{2})^{m^0} + (\frac{3}{2})^{m^0} - 1)*3+1}{2^{e^0}}$$

$$k^1*2^{e^0+m^0+1} + 2^{e^0+m^0} = (k^0+1)*3^{m^0+1} - 2^{m^0}$$

$$(2*k^1+1)*2^{e^0} = 2*f*3^{m^0+1} + 3^{m^0+1} - 1 \quad f = 0,1,2,3,\cdots\cdots$$

$$k^1 = k^0 * \frac{3^{m^0+1}}{2^{e^0+m^0+1}} + \frac{3^{m^0+1} - 2^{m^0} - 2^{e^0+m^0}}{2^{e^0+m^0+1}} \tag{9}$$

So in the main cycle we can get three number sequence $k^n, m^n, e^n, (n = 0,1,2,\cdots\cdots)$.

$$k^n = k^0 * \prod_{i=0}^{n-1} \frac{3^{m^i+1}}{2^{e^i+m^i+1}} + \sum_{j=0}^{n-2} \frac{3^{m^j+1} - 2^{m^j} - 2^{e^j+m^j}}{2^{e^j+m^j+1}} \prod_{i=j+1}^{n-1} \frac{3^{m^i+1}}{2^{e^i+m^i+1}} + \frac{3^{m^{n-1}+1} - 2^{m^{n-1}} - 2^{e^{n-1}+m^{n-1}}}{2^{e^{n-1}+m^{n-1}+1}} \tag{10}$$

If not the only cycle when it get n.

$$k' = k^{\alpha*n} \quad \alpha = 0,1,2,3$$

$$n = 1$$

$$k' = k^1 = k^0 * \frac{3^{m^0+1}}{2^{e^0+m^0+1}} + \frac{3^{m^0+1} - 2^{m^0} - 2^{e^0+m^0}}{2^{e^0+m^0+1}}$$

$n > 1$

So as we know

$$k' = k^n = k^0 * \prod_{i=0}^{n-1} \frac{3^{m^i+1}}{2^{e^i+m^i+1}} + \sum_{j=0}^{n-2} \frac{3^{m^j+1} - 2^{m^j} - 2^{e^j+m^j}}{2^{e^j+m^j+1}} \prod_{i=j+1}^{n-1} \frac{3^{m^i+1}}{2^{e^i+m^i+1}} + \frac{3^{m^{n-1}+1} - 2^{m^{n-1}} - 2^{e^{n-1}+m^{n-1}}}{2^{e^{n-1}+m^{n-1}+1}}$$

$$k'(2^{e^0+m^0+1} - 3^{m^0+1}) = 3^{m^0+1} - 2^{m^0} - 2^{e^0+m^0} \tag{11}$$

When $k'$ is even, so $m^0 = 0$ and $e^0 > 0$, then we can find the left is below zero, but $2^{e^0+m^0+1} - 3^{m^0+1} > 0$. So $k' < 0$ this is not obey. So that all $k' = 0$.

When $k'$ is odd, so $m^0 > 0$ , then we can find the left is below zero, we can only get $k' = -1, -3$ this two integer. But this is not obey. So that all $k'$ is not odd.

So we can get if n=1 $k' = 0$, and problem is OK.

$$k^n * \prod_{i=0}^{n-1} 2^{e^i+m^i+1} = k^0 * \prod_{i=0}^{n-1} 3^{m^i+1} + \prod_{i=0}^{n-1} 2^{e^i+m^i+1} * \left(\sum_{j=0}^{n-2} \frac{3^{m^j+1} - 2^{m^j} - 2^{e^j+m^j}}{2^{e^j+m^j+1}} \prod_{i=j+1}^{n-1} \frac{3^{m^i+1}}{2^{e^i+m^i+1}} + \frac{3^{m^{n-1}+1} - 2^{m^{n-1}} - 2^{e^{n-1}+m^{n-1}}}{2^{e^{n-1}+m^{n-1}+1}}\right)$$

$$k' * \left(\prod_{i=0}^{n-1} 2^{e^i+m^i+1} - \prod_{i=0}^{n-1} 3^{m^i+1}\right) = \prod_{i=0}^{n-2} 2^{e^i+m^i+1} * (3^{m^{n-1}+1} - 2^{m^{n-1}} - 2^{e^{n-1}+m^{n-1}}) + \sum_{j=1}^{n-2} \left((3^{m^j+1} - 2^{m^j} - 2^{e^j+m^j}) \prod_{i=j+1}^{n-1} 3^{m^i+1} \prod_{i=0}^{j-1} 2^{e^i+m^i+1}\right) + ((3^{m^0+1} - 2^{m^0} - 2^{e^0+m^0}) \prod_{i=1}^{n-1} 3^{m^i+1} \tag{12}$$

If all $k'$ is even, so $m^0 = 0$. When $e^0 > 1$, then we can find the right is below zero, but $\prod_{i=0}^{n-1} 2^{e^i+m^i+1} - \prod_{i=0}^{n-1} 3^{m^i+1} > 0$. So $k' < 0$ this is not obey; when $e^0 = 1$, the right is zero and $\prod_{i=0}^{n-1} 2^{e^i+m^i+1} - \prod_{i=0}^{n-1} 3^{m^i+1} \neq 0$, so all $k' = 0$.

So if there have one or more $k'$ is odd, the problem can't found. In other world, the $\eta$ kind number is not zero.

## 6 The number of $\eta$ between The number of $\gamma$

Because when we use the $z = x + s - 1$, we can get a real simple function for $T(z)$.

$$\mathrm{T}(z) = \frac{(2*(1-\mathrm{x}+\mathrm{s})+1)*(\mathrm{x}+\mathrm{s}-1)+1-(x-s)}{2} \tag{13}$$

As $z$ is an odd number, and $x_{odd} = s_{odd}$.

$$\mathrm{T}(z_{odd}) = \frac{3*(\mathrm{x}_{odd}+\mathrm{s}_{odd}-1)+1}{2}$$

As $z$ is an even number, and $x_{even} = s_{even} + 1$

$$\mathrm{T}(z_{even}) = \frac{\mathrm{x}_{even}+\mathrm{s}_{even}-1}{2}$$

And when $T(z_0) = z_1$

$$x_1 + s_1 - 1 = \frac{(2*(1-\mathrm{x}_0+\mathrm{s}_0)+1)*(\mathrm{x}_0+\mathrm{s}_0-1)+1-(x_0-s_0)}{2}$$

Simplified it

$$x_1 + s_1 = (x_0 + s_0) + x_0 - x_0^2 + s_0^2 \tag{14}$$

So we can get $T^n(z_0) = z_n$

$$x_n + s_n = (x_0 + s_0) + \sum_{j=0}^{n-1} x_j - \sum_{j=0}^{n-1} x_0^2 + \sum_{j=0}^{n-1} s_0^2$$

Simplified it

$$x_n + s_n = (x_0 + s_0) + \sum_{j=0}^{n-1} x_j + \sum_{j=0}^{n-1} (s_j + x_j)(s_j - x_j)$$

If we set that $z_n = z_0$ ,we can get that $x_n = x_0$ , $s_n = s_0$.

$$0 = \sum_{j=0}^{n-1} x_j + \sum_{j=0}^{n-1} (s_j + x_j)(s_j - x_j) \tag{15}$$

If we summed that there a cycle don't contain $\eta$ kind number, and the cycle contain only one

$\gamma$ number in one time cycle.

Base on the rule $T(z)$, we can get the cycle as $\alpha \to \beta \to \alpha$

We can set $z_0 = \alpha$, $z_1 = \beta$.

Base on $T(z)$ we can get $z_0 \to z_1$

$$s_0 = x_0 = \frac{s_1 + x_1}{3}$$

$$0 = -\frac{x_0 - 1}{2} + \sum_{j=2}^{n-1} x_j + \sum_{j=2}^{n-1} (s_j + x_j)(s_j - x_j) \tag{16}$$

And we can find when and only when all $x_i = 1(i = 0,1,2\cdots)$, the cycle can be right.

If the cycle add a $\gamma$ kind number, there must get a part of the cycle is $\alpha \to \gamma \to \beta \to \alpha$.

We can set $z_{n-1} = \gamma$, $z_0 = \beta$, $z_1 = \alpha$.

Base on $T(z)$ we can get $z_{n-1} \to z_0 \to z_1$

$s_{n-1} + 1 = x_{n-1} = s_0 + x_0$ and $s_{n-1} = 2s_0$

$s_0 + 1 = x_0 = s_1 + x_1$ and $s_1 = x_1$

So $s_0 = 2 * x_1 - 2$

$$0 = -5x_1 + 6 + \sum_{j=2}^{n-2} x_j + \sum_{j=2}^{n-2} (s_j + x_j)(s_j - x_j) \tag{17}$$

and $x_1$ is an even number

So that we can find if the cycle add one or more $\gamma$ kind number, the cycle can't be get, because the right side must be a negative number.

So in the contrast way, the cycle add a $\eta$ kind number, there must get a part of the cycle is $\alpha \to \beta \to \eta \to \alpha$

We can set $z_{n-1} = \beta$, $z_0 = \eta$, $z_1 = \alpha$.

Base on $T(z)$ we can get $z_{n-1} \to z_0 \to z_1$

$s_{n-1} + 1 = x_{n-1} = s_n + x_n = 2 * x_0$

$$s_0 = x_0 = \frac{2}{3} s_1 = \frac{2}{3} x_1$$

$$0 = \frac{x_1}{3} + 1 + \sum_{j=2}^{n-2} x_j + \sum_{j=2}^{n-2} (s_j + x_j)(s_j - x_j) \qquad (18)$$

So that we can find if the cycle add one or more $\eta$ kind number too, the cycle can't be get, because the right side must be a positive number.

Overall, we can know base on this cycle, if the number of $\eta$ kind number is zero the 3x+1 will be proofed.

## ACKNOWLEDGMENTS

Thanks for my middle school teacher, he told this problem when I was a middle school student. So I get a great haptic in science to find something.

## REFERENCES

[1] Peter Eisele and Karl-Peter Hadeler, Game of Cards, Dynamical Systems, and a Characterization of the Floor and Ceiling Functions, Amer. Math. Monthly 97 (1990), 466–477. (MR 91h: 58086).

[2] Gaston Gonnet, Computations on the 3n+1 Conjecture, MAPLE Technical Newsletter 0, No. 6, Fall 1991.

[3] Edward Belaga and Maurice Mignotte, The Collatz problem and its generalizations: Experimental Data. Table 1. Primitive cycles of 3x + d mappings, Univ. of Strasbourg preprint 2006-015, 9 pages+400+ page table.